\documentclass[a4paper,11pt]{amsart}
\usepackage{ctex}
\usepackage{amssymb}
\usepackage{latexsym}
\usepackage{exscale}
\usepackage{amsfonts}
\usepackage{graphicx}
\usepackage{mathrsfs}
\usepackage{amsmath,amscd,amsthm}
\usepackage{bbm,pifont}
\usepackage{enumerate}
\usepackage[colorlinks, citecolor=blue,  hypertexnames=false]{hyperref}
\usepackage{amssymb}
\usepackage{latexsym}
\usepackage{exscale}

\allowdisplaybreaks

\scrollmode
\allowdisplaybreaks
\numberwithin{equation}{section}
\newtheorem{theorem}{Theorem}[section]

\newtheorem{definition}[theorem]{Definition}

\allowdisplaybreaks

\begin{document}

\title[Flag Lipschitz spaces on $\mathbb{H}^{n}$]{
Marcinkiewicz multipliers and Lipschitz spaces on Heisenberg groups}

\author{Yanchang Han}
\address{School  of Mathematic Sciences,
South China Normal University,
Guangzhou, 510631, P.R. China}
\email{20051017@m.scnu.edu.cn}

\author{Yongsheng Han}
\address{Department of Mathematics,
Auburn University,
Auburn, AL 36849, U.S.A.}
\email{hanyong@mail.auburn.edu}

\author{Ji Li}
\address{Department of Mathematics,  Macquarie University,  NSW, 2109, Australia.}
\email{ji.li@mq.edu.au}

\author{Chaoqiang Tan}
\address{Department of Mathematics, Shantou University, Shantou, Guangdong 515063, China}
\email{cqtan@stu.edu.cn}

\thanks{The first  author is supported by National Natural Science Foundation of China  (Grant No. 11471338) and
Guangdong Province Natural Science Foundation  (Grant No. 2014A030313417); The third author is supported by the
Australian Research Council under Grant No.~ARC-DP160100153 and by Macquarie University New Staff Grant}

\maketitle

\begin{abstract}
The Marcinkiewicz multipliers are $L^{p}$ bounded for $1<p<\infty $ on the
Heisenberg group $\mathbb{H}^{n}\simeq \mathbb{C}^{n}\times \mathbb{R}$ (M\"uller, Ricci and Stein \cite{MRS}). This is
surprising in the sense that these multipliers are invariant under a two parameter
group of dilations on $\mathbb{C}^{n}\times \mathbb{R}$, while there is
\emph{no} two parameter group of \emph{automorphic} dilations on $\mathbb{H}
^{n}$. The purpose of this paper is to establish a theory of the flag Lipschitz space on the Heisenberg group
$\mathbb{H}^{n}\simeq \mathbb{C}^{n}\times \mathbb{R}$ in the sense
`intermediate' between the classical Lipschitz space on the Heisenberg group
$\mathbb{H}^{n}$ and the product Lipschitz space on
$\mathbb{C}^{n}\times \mathbb{R}$. We characterize this flag Lipschitz space
via the Littelewood-Paley theory and prove
that flag singular integral operators, which include the
Marcinkiewicz multipliers, are bounded on these flag Lipschitz spaces.
\end{abstract}
\maketitle

\section{Introduction}

Classical Calder\'{o}n-Zygmund singular integrals
commute with
the one parameter dilations on $\mathbb{R}^{n}$, $\delta \cdot x=(\delta
x_{1},...,\delta x_{n})$ for $\delta >0$, while the \emph{product}
Calder\'{o}n-Zygmund singular integrals commute with the
multi-parameter dilations on $\mathbb{R}^{n}$, $\delta \cdot x=(\delta
_{1}x_{1},...,\delta _{n}x_{n})$ for $\delta =(\delta _{1},...,\delta
_{n})\in \mathbb{R}_{+}^{n}$.

In the product Calder\'{o}n-Zygmund theory, the product singular integral
operators are of the form $Tf=K\ast f,$ where $K$ is homogeneous, that is, ${
\delta _{1}}...{\delta _{n}}K(\delta \cdot x)=K(x),$ or, more generally, $
K(x)$ and ${\delta _{1}}...{\delta _{n}}K(\delta
\cdot x)$ satisfy the same size, smoothness and cancellation
conditions. Such operators have been studied for
example in Gundy-Stein \cite{GS}, R. Fefferman and Stein \cite{FS}, R.
Fefferman (\cite{F1, F2, F3}), Chang \cite{Cha}, Chang and R. Fefferman
(\cite{CF1}, \cite{CF2}, \cite{CF3}), Journ\'{e} (\cite{J1}, \cite{J2}), and
Pipher \cite{P}. More precisely, R. Fefferman and Stein \cite{FS} studied the $
L^{p}$ boundedness ($1<p<\infty $) for the product convolution singular
integral operators. Journ\'{e} in (\cite{J1, J3}) introduced a non-convolution
product singular integral operators and established the product $T1$ theorem
and proved the $L^{\infty }\rightarrow BMO$ boundedness for such operators.
The product Hardy space $H^{p}\left( \mathbb{R}^{n}\times \mathbb{R}
^{m}\right) $ was first introduced by Gundy and Stein \cite{GS}. Chang and
R. Fefferman (\cite{CF1, CF2, CF3}) developed the theory of atomic
decomposition and established the dual space of Hardy space $H^{1}\left(
\mathbb{R}^{n}\times \mathbb{R}^{m}\right) $, namely the product $BMO\left(
\mathbb{R}^{n}\times \mathbb{R}^{m}\right) $ space. Carleson \cite{Car} disproved
by a counter-example a conjecture that the product atomic Hardy space on $
\mathbb{R}^{n}\times \mathbb{R}^{m}$ could be defined by rectangle atoms.
This motivated Chang and R. Fefferman to replace the role of cubes in the
classical atomic decomposition of $H^{p}\left( \mathbb{R}^{n}\right) $ by
arbitrary open sets of finite measures in the product $H^{p}\left( \mathbb{R}
^{n}\times \mathbb{R}^{m}\right) $. Subsequently, R. Fefferman in \cite{F2}
established the criterion of the $H^{p}\rightarrow L^{p}$ boundedness of
singular integral operators in Journ\'{e}'s class by considering its action
only on rectangle atoms via Journ\'{e} lemma. However, R. Fefferman's
criterion cannot be extended to three or more parameters without further
assumptions on the nature of $T$ as shown in Journ\'{e} (\cite{J1, J3}). In
fact, Journ\'{e} provided a counter-example in the three-parameter setting
of singular integral operators such that R. Fefferman's criterion breaks
down. Subsequently, the $H^{p}$ to $L^{p}$ boundedness for Journ\'{e}'s
class of singular integral operators with arbitrary number of parameters was
established by J. Pipher \cite{P} by considering directly the action of the
operator on (non-rectangle) atoms and an extension of Journ\'{e}'s geometric
lemma to higher dimensions.

On the other hand, multi-parameter analysis has only recently been developed
for $L^{p}$ theory with $1<p<\infty $ when the underlying multi-parameter
structure is not explicit, but \emph{implicit}, as in the flag
multi-parameter structure studied on the
Heisenberg group $\mathbb{H}^{n}$ by M\"uller, Ricci and Stein in (\cite{MRS},\cite{MRS2}).
See also Phong and Stein in \cite{PS} and Nagel, Ricci and Stein in \cite{NRS}. In \cite{MRS} and \cite{MRS2},
the authors obtained the surprising result that
certain Marcinkiewicz multipliers, invariant under a two-parameter group of
dilations on $\mathbb{C}^{n}\times \mathbb{R}$, are bounded on $L^{p}\left(
\mathbb{H}^{n}\right) $, \emph{despite} the absence of a two-parameter
automorphic group of dilations on $\mathbb{H}^{n}$. To be precise,
the Heisenberg group $\mathbb{H}^{n}$ is the one consists of the set
$$\mathbb{C}^{n} \times \mathbb{R}= \{[z,t]: z \in \mathbb{C}^{n}, t \in \mathbb{R}\}$$
with the multiplication law

$$[z,t]\circ [z',t']=\Big[z+z',t+t'+2 Im \big({z}{\bar{z'}}\big) \Big],$$
where identity is the origin [0,0] and the inverse is given by $[z, t]^{-1} = [-z,-t]$.

In addition to the Heisenberg group multiplication law, nonisotropic dilations of $\mathbb{H}^{n}$ are given by
\begin{equation*}
\delta_{r} : \mathbb{H}^{n} \rightarrow \mathbb{H}^{n}, \quad \delta_{r} ([z,t]) = [r z, r^{2}t].
\end{equation*}
A trivial computation shows that $\delta_{r}$ is an automorphism of $\mathbb{H}^{n}$ for every $r>0$. However, the standard isotropic dilations of $\mathbb R^{2n+1}$ are not automorphisms of $\mathbb{H}^{n}.$

The ``norm'' function $\rho$ on $\mathbb{H}^{n}$ is defined by
\begin{equation*}
\rho([z,t]):= (|z|^2+ |t|)^{1/2}.
\end{equation*}
It is easy to see that $\rho([z,t]^{-1}) = \rho([-z,-t])  =\rho([z,t])$,  $\rho(\delta_{r}([z,t])) = r\rho([z,t])$,
$\rho([z,t]) = 0$ if and only if $ [z,t]=[0,0]$,  and
\begin{equation*}
\rho([z,t]\circ [z',t'])  \leq \gamma (\rho([z,t]) + \rho([z',t'])),
\end{equation*}
where $\gamma >1$ is a constant.

The Haar measure on $\mathbb{H}^{n}$ is known to just coincide with the Lebesgue measure on $\mathbb{R}^{2n+1}$.
The vector fields
\begin{align*}
T: =\frac{\partial}{\partial t}, \quad&X_{j}:=\frac{\partial}{\partial x_{j}} -2y_{j}\frac{\partial}{\partial t}, \quad
Y_{j}:=\frac{\partial}{\partial y_{j}} +2x_{j}\frac{\partial}{\partial t}, \quad j =1,\cdots, n,
\end{align*}
form a natural basis for the Lie algebra of left-invariant vector fields on $\mathbb{H}^{n}$.
The standard sub-Laplacian $\mathcal{L}$ on
the Heisenberg group is defined by
\begin{equation*}
\mathcal{L}\equiv -\sum_{j=1}^{n}\left( X_{j}^{2}+Y_{j}^{2}\right).
\end{equation*}
The operators $\mathcal{L}$ and $T=\frac{\partial }{\partial t}$ commute,
and so do their spectral measures $dE_{1}\left( \xi \right) $ and $
dE_{2}\left( \eta \right) $. Given a bounded function $m\left( \xi ,\eta
\right) $ on $\mathbb{R}_{+}\times \mathbb{R}$, define the multiplier
operator $m\left( \mathcal{L},iT\right) $ on $L^{2}\left( \mathbb{H}
^{n}\right) $ by
\begin{equation*}
m\left( \mathcal{L},iT\right) =\int_{\mathbb R}
\int_{\mathbb{R}_{+}} m\left( \xi ,\eta \right) dE_{1}\left( \xi
\right) dE_{2}\left( \eta \right) .
\end{equation*}

 Then $m\left( \mathcal{L},iT\right) $ is automatically bounded on
$ L^{2}\left( \mathbb{H}^{n}\right) $, and if one imposes
Marcinkiewicz conditions on the multiplier $m\left( \xi ,\eta
\right)$, namely
\begin{equation*}
\left\vert \left( \xi \partial _{\xi }\right) ^{\alpha }\left( \eta \partial
_{\eta }\right) ^{\beta }m\left( \xi ,\eta \right) \right\vert \leq
C_{\alpha ,\beta },
\end{equation*}
for all $\left\vert \alpha \right\vert ,\beta \geq 0$,
then M\"uller, Ricci and Stein \cite{MRS} proved that the Marcinkiewicz multiplier operator
$m\left( \mathcal{L}, iT\right) $ is a bounded operator on $L^{p}\left( \mathbb{H}
^{n}\right) $ for $1<p<\infty$.
This is surprising since these multipliers are invariant under a two parameter
group of dilations on $\mathbb{C}^{n}\times \mathbb{R}$, while there is
\emph{no} two parameter group of \emph{automorphic} dilations on $\mathbb{H}
^{n}$. Moreover, they showed that Marcinkiewicz multiplier can be characterized by convolution operator with the form $f\ast K$ where, however, $K$ is a flag convolution kernel. A flag convolution kernel on $\mathbb{H}^{n}=\mathbb{C}
^{n}\times \mathbb{R}$ is a distribution $K$ on $\mathbb{H}^{n}$ which
coincides with a $C^{\infty }$ function away from the coordinate subspace $
\{(0,u)\}\subset \mathbb{H}^{n}$, where $0\in \mathbb{C}^{n}$ and $u\in
\mathbb{R}$, and satisfies

\begin{enumerate}
\item (Differential Inequalities) For any multi-indices $\alpha =(\alpha
_{1},\cdots ,\alpha _{n})$,
\begin{equation*}
\left\vert \partial _{z}^{\alpha }\partial _{u}^{\beta }K(z,u)\right\vert
\leq C_{\alpha ,\beta }\left\vert z\right\vert ^{-2n-|\alpha |}\cdot \left(
\left\vert z\right\vert ^{2}+\left\vert u\right\vert \right) ^{-1-\beta }
\end{equation*}
for all $(z,u)\in \mathbb{H}^{n}$ with $z\neq 0$ and all $|\alpha|, \beta\geq 0.$

\item (Cancellation Condition)
\begin{equation*}
\left\vert \int_{\mathbb{R}}\partial _{z}^{\alpha }K(z,u)\phi _{1}(\delta
u)du\right\vert \leq C_{\alpha }|z|^{-2n-|\alpha |}
\end{equation*}
for every multi-index $\alpha $ and every normalized bump function $\phi
_{1} $ on $\mathbb{R}$ and every $\delta >0$;
\begin{equation*}
\left\vert \int_{\mathbb{C}^{n}}\partial _{u}^{\beta }K(z,u)\phi _{2}(\delta
z)dz\right\vert \leq C_{\gamma }|u|^{-1-\beta }
\end{equation*}
for every index $\beta\geq 0 $ and every normalized bump function $\phi _{2}$
on $\mathbb{C}^{n}$ and every $\delta >0$; and
\begin{equation*}
\left\vert \int_{\mathbb{H}^{n}}K(z,u)\phi _{3}(\delta _{1}z,\delta
_{2}u)dzdu\right\vert \leq C
\end{equation*}
for every normalized bump function $\phi _{3}$ on $\mathbb{H}^{n}$ and every
$\delta _{1}>0$ and $\delta _{2}>0$.
\end{enumerate}
They also proved that flag singular integral operators are bounded
on $L^{p}\left( \mathbb{H}^{n}\right) $ for $1<p<\infty$. See \cite{NRSW1} and \cite{NRSW2} for more about flag singular integral operators on homogeneous groups.

At the endpoint estimates, it is natural to expect that the boundedness on Hardy and $\mathrm{BMO}$ spaces are available. However, the lack of automorphic dilations underlies the failure of such
multipliers to be in general bounded on the classical Hardy space $H^{1}$ and also precludes a pure product Hardy space theory on the Heisenberg group.
This was the original motivation in \cite{HLS} to develop a theory of \emph{flag} Hardy
spaces $H_{flag}^{p}$ on the Heisenberg group, $0<p\leq 1$, that is, in a
sense `intermediate' between the classical Hardy spaces $H^{p}(\mathbb{H}
^{n})$ and the
product Hardy spaces $H_{product}^{p}(\mathbb{C}^{n}\times \mathbb{R})$
(A. Chang and R. Fefferman (\cite{CF1}, \cite{CF2}, \cite{F1}, \cite{F2},
\cite{F3})). They showed that singular integrals with flag kernels, which include the
aforementioned Marcinkiewicz multipliers, are bounded on $H_{flag}^{p}$, as
well as from $H_{flag}^{p}$ to $L^{p}$, for $0<p\leq 1$. Moreover, they constructed a singular integral with a flag kernel on the Heisenberg group, which is not bounded on the classical Hardy spaces $H^{1}(\mathbb{H}
^{n}).$ Since, as pointed out in \cite{HLS}, the flag Hardy space $H_{flag}^{p}(\mathbb{H}
^{n})$ is contained in the classical Hardy space $H^{p}(\mathbb{H}
^{n}),$ this counterexample implies that $H_{flag}^{1}(\mathbb{H}
^{n})\subsetneqq H^{1}(\mathbb{H}^{n}).$

A basic question arises: Can one establish the endpoint estimates of flag singular integral operators on the Lipschitz space?
It was well known that the Lipschitz spaces are important for the study of partial differential equations. These spaces are also crucial in developing  the so-called first order calculus, systematic theories on metric measure spaces since the end of the 1970s. See, for example, \cite{C, He1, He2} and for more about the classical Lipschitz spaces see \cite{HSV}, \cite{JTW}, \cite{K1}, \cite{K2}, \cite{MR}, \cite{S1} and \cite{S2}.

The purpose of this paper is to answer this question. More precisely, we will establish a theory of the flag Lipschitz space on Heisenberg group
$\mathbb{H}^{n}\simeq \mathbb{C}^{n}\times \mathbb{R}$ in the sense
`intermediate' between the classical Lipschitz space on Heisenberg group
$\mathbb{H}^{n}$ and the product Lipschitz space on
$\mathbb{C}^{n}\times \mathbb{R}$.

We will characterize the flag Lipschitz space
via the Littelewood-Paley theory and prove
that flag singular integral operators, which include the
Marcinkiewicz multipliers, are bounded on the flag Lipschitz space on Heisenberg group
$\mathbb{H}^{n}$.

Now we introduce the following notations:
$$\Delta^1_{(u,v)}(f)(z,r)=f((z,r)\circ(u,v)^{-1})-f(z,r),$$
$$\Delta^{1,Z}_{(u,v)}(f)(z,r)=f((z,r)\circ(u,v))+f((z,r)\circ(u,v)^{-1})-2f(z,r),$$
and
$$\Delta^2_{w}(f)(z,r)=f(z,r-w)-f(z,r),$$
$$ \Delta^{2,Z}_{w}(f)(z,r)=f(z,r+w)+f(z,r-w)-2f(z,r).$$

The flag Lipschitz space on Heisenberg group is defined by the following
\begin{definition}\label{lip} A continuous function $f(z,r)$ defined on $\mathbb H^n$ belongs to the flag Lipschitz space $Lip_{flag}^\alpha$ with $\alpha=(\alpha_1,\alpha_2), \alpha_1, \alpha_2 >0$ if and only if

(i) when $0<\alpha_1, \alpha_2<1,$
\begin{eqnarray}
\label{lip1}&&|\Delta^2_{w}\Delta^1_{(u,v)}(f)(z,r)|\\
&=&|f((z,r)\circ(u,v+w)^{-1})-f((z,r)\circ(u,v)^{-1})-f((z,r)\circ(0,w)^{-1})+f(z,r)|\nonumber\\
&\leq& C|(u,v)|^{\alpha_1}|w|^{\alpha_2},\nonumber\end{eqnarray}
where $|(u,v)|^2=|u|^2+|v|;$

(ii) when $\alpha_1=1, 0<\alpha_2<1,$
\begin{eqnarray}
\label{lip2}&&|\Delta^2_{w}\Delta^{1,Z}_{(u,v)}(f)(z,r)|\\
&=&|[f((z,r)\circ(u,v))+f((z,r)\circ(u,v)^{-1})-2f(z,r)]\nonumber\\
&&-[f((z,r)\circ(u,v-w))+f((z,r)\circ(u,v+w)^{-1})-2f((z,r)\circ(0,w)^{-1})]|\nonumber\\
&\leq& C|(u,v)||w|^{\alpha_2};\nonumber \end{eqnarray}

(iii) when $0<\alpha_1<1, \alpha_2=1,$
\begin{eqnarray}
\label{lip3}&&|\Delta^{2,Z}_{w}\Delta^1_{(u,v)}(f)(z,r)|\\
&=&|[f((z,r)\circ(0,w))+f((z,r)\circ(0,w)^{-1})-2f(z,r)]\nonumber\\
&&-[f((z,r)\circ(u,-v+w))+f((z,r)\circ(u,v+w)^{-1})-2f((z,r)\circ(u,v)^{-1})]|\nonumber\\
&\leq& C|(u,v)|^{\alpha_1}|w|;\nonumber\end{eqnarray}

(iv) when $\alpha_1=\alpha_2=1,$
\begin{eqnarray}
\label{lip4}&&|\Delta^{2,Z}_{w}\Delta^{1,Z}_{(u,v)}(f)(z,r)|\\
&=&|[f((z,r)\circ(u,v+w))+f((z,r)\circ(u,v-w)^{-1})-2f((z,r)\circ(0,w))]\nonumber\\
&&+[f((z,r)\circ(u,v-w))+f((z,r)\circ(u,v+w)^{-1})-2f((z,r)\circ(0,w)^{-1})]\nonumber\\
&&-2[f((z,r)\circ(u,v))+f((z,r)\circ(u,v)^{-1})-2f(z,r))]|\nonumber\\
&\leq&C|(u,v)| |w|;\nonumber\end{eqnarray} for all $(z,r)$ and
$(u,v) \in \mathbb H^n, w\in \mathbb R$ and the constant $C$ is
independent of $(z,r), (u,v)$ and $w.$

When $\alpha=(\alpha_1, \alpha_2)$ with $\alpha_1, \alpha_2>1,$ we write $\alpha_1=m_1+r_1$ and $\alpha_2=m_2+r_2$ where $m_1, m_2$ are integers and $0<r_1, r_2\leq 1. f(z,u)\in Lip_{flag}^\alpha$ means that $f(z,u)$ is a $C^{m_1+m_2}$ function, modulo polynomials of degree not exceeding $m_1+m_2,$ such that all the partial derivatives ${\partial_z^{m_1}\partial_u^{m_2} f(z,u)}$ belong to $Lip_{flag}^r$ with $r=(r_1,r_2).$

If $f\in Lip_{flag}^\alpha, ||f||_{Lip_{flag}^\alpha},$ the norm of $f,$ is defined as the smallest constant $C$ such that \eqref{lip1} to \eqref{lip4} hold.
\end{definition}

Note that the flag structure is involved in the the definition of $Lip_{flag}^{\alpha}(\mathbb H^n)$ with $\alpha=(\alpha_1,\alpha_2)$ and Zygmund conditions are used when $\alpha_1=1$ or $\alpha_2=1,$ or both $\alpha_1=\alpha_2=1.$

In order to obtain the boundedness of flag singular integral operators on the flag Lipschitz space $Lip_{flag}^\alpha$, we characterize
$Lip_{flag}^\alpha$ via the Littlewood-Paley theory.
We begin with recalling the standard Calder\'{o}n reproducing formula on the Heisenberg
group. Note that spectral theory was used in place of the Calder\'{o}n
reproducing formula in \cite{MRS2}.

\begin{theorem}
\label{GM}(Corollary 1 of \cite{GM}) There is a radial function $\psi \in C^{\infty }\left(
\mathbb{H}^{n}\right) $ satisfying $\psi \in \mathcal{S}\left( \mathbb{H}^{n}\right)$
and all moments of $\psi$ vanish,
such that the following Calder\'{o}n reproducing formula holds:
\begin{equation*}
f=\int_{0}^{\infty }\psi _{s}\ast \psi _{s}\ast f\frac{ds}{s},\ \ \
\ \ f\in L^{2}\left( \mathbb{H}^{n}\right) ,
\end{equation*}
where $\ast $ is Heisenberg convolution and $\psi _{s}\left( z,t\right)
=s^{-2n-2}\psi \left( \frac{z}{s},\frac{t}{s^{2}}\right) $ for $s>0$.
\end{theorem}

We now wish to extend this formula in two ways: (1) encompass the flag structure on the
Heisenberg group $\mathbb{H}^{n}$; (2) converges in the distribution space. For this purpose,
following \cite{MRS2}, we construct a Littlewood-Paley \emph{component}
function $\psi $ defined on $\mathbb{H}^{n}\simeq \mathbb{C}^{n}\times
\mathbb{R}$, given by the partial convolution $\ast _{2}$ in the second
variable only:
\begin{equation*}
\psi (z,u)=\psi ^{(1)}\ast _{2}\psi ^{(2)}(z,u)=\int_{\mathbb{R}}\psi
^{(1)}(z,u-v)\psi ^{(2)}(v)dv,\ \ \ \ \ \left( z,u\right) \in \mathbb{C}
^{n}\times \mathbb{R},
\end{equation*}
where $\psi ^{(1)}\in \mathcal{S}(\mathbb{H}^{n})$ is as in Theorem
\ref{GM}, and $\psi ^{\left( 2\right) },$ a real even function defined in $\mathcal{S}\left( \mathbb{R}
\right) $ satisfies
\begin{equation*}
\int_{0}^{\infty }|\widehat{{\psi }^{(2)}}(t\eta )|^{2}\frac{dt}{t}=1
\end{equation*}
for all $\eta \in \mathbb{R}\backslash \{0\}$, along with the moment
conditions
\begin{eqnarray*}
\int\limits_{\mathbb{H}^{n}}z^{\alpha }u^{\beta }\psi ^{(1)}(z,u)dzdu &=&0,\
\ \ \ \ \left\vert \alpha \right\vert, \beta \geq 0, \\
\int\limits_{\mathbb{R}}v^{\gamma }\psi ^{(2)}(v)dv &=&0,\ \ \ \ \ \gamma
\geq 0.
\end{eqnarray*}
Thus we have
\begin{equation}
f(z,u)=\int_{0}^{\infty }\int_{0}^{\infty }{\psi}_{s,t}\ast \psi
_{s,t}\ast f(z,u)\frac{ds}{s}\frac{dt}{t},\ \ \ \ \ f\in L^{2}\left( \mathbb{
H}^{n}\right) ,  \label{CRF st}
\end{equation}
where the functions $\psi _{s,t}$ are given by
\begin{equation}
\psi _{s,t}\left( z,u\right) =\psi _{s}^{(1)}\ast _{2}\psi _{t}^{(2)}\left(
z,u\right) ,  \label{defpsits}
\end{equation}
with
\begin{equation*}
\psi _{s}^{(1)}(z,u)=s^{-2n-2}\psi ^{(1)}(\frac{z}{s},\frac{u}{s^{2}})\,\,
\text{and}\,\,\psi _{t}^{(2)}(v)=t^{-1}\psi ^{(2)}(\frac{v}{t}),
\end{equation*}
and where the integrals in (\ref{CRF st}) converge in $L^{2}\left( \mathbb{H}
^{n}\right) $. Indeed,
\begin{eqnarray*}
{\psi}_{s,t}\ast _{\mathbb{H}^{n}}\psi _{s,t}\ast _{\mathbb{H}
^{n}}f\left( z,u\right) &=&\left( {\psi}_{s}^{(1)}\ast _{2}\psi
_{t}^{(2)}\right) \ast _{\mathbb{H}^{n}}\left( {\psi}_{s}^{(1)}\ast
_{2}\psi _{t}^{(2)}\right) \ast _{\mathbb{H}^{n}}f\left( z,u\right) \\
&=&\left( {\psi}_{s}^{(1)}\ast _{\mathbb{H}^{n}}{\psi}
_{s}^{(1)}\right) \ast _{\mathbb{H}^{n}}\left( \left(\psi _{t}^{(2)}\ast _{\mathbb{
R}}\psi _{t}^{(2)}\right) \ast _{2}f\left( z,u\right)\right)
\end{eqnarray*}
yields (\ref{CRF st}) upon invoking the standard Calder\'{o}n reproducing
formula on $\mathbb{R}$ and then Theorem \ref{GM} on $\mathbb{H}^{n}$:
\begin{eqnarray*}
&&\int_{0}^{\infty }\int_{0}^{\infty }{\psi}_{s,t}\ast _{\mathbb{H}
^{n}}\psi _{s,t}\ast _{\mathbb{H}^{n}}f(z,u)\frac{ds}{s}\frac{dt}{t} \\
&=&\int_{0}^{\infty }{\psi}_{s}^{(1)}\ast _{\mathbb{H}^{n}}{\psi}
_{s}^{(1)}\ast _{\mathbb{H}^{n}}\left\{ \int_{0}^{\infty }\psi
_{t}^{(2)}\ast _{\mathbb{R}}\psi _{t}^{(2)}\ast _{2}f\left( z,u\right) \frac{
dt}{t}\right\} \frac{ds}{s} \\
&=&\int_{0}^{\infty }{\psi}_{s}^{(1)}\ast _{\mathbb{H}^{n}}{\psi}
_{s}^{(1)}\ast _{\mathbb{H}^{n}}f\left( z,u\right) \frac{ds}{s}=f\left(
z,u\right) .
\end{eqnarray*}
Note that the Littlewood-Paley component
function $\psi $ satisfies the \emph{flag} moment conditions, so-called because they
include only half of the product moment conditions associated with the
product $\mathbb{C}^{n}\times \mathbb{R}$:
\begin{equation}
\int\limits_{\mathbb{R}}u^{\alpha }\psi (z,u)du=0,\ \ \ \ \ \text{for all }
\alpha \geq 0 \text{ and }z\in \mathbb{C}^{n}.  \label{momcond}
\end{equation}
Indeed, with the change of variable $u^{\prime }=u-v$ and the binomial
theorem
\begin{equation*}
\left( u^{\prime }+v\right) ^{\beta }=\sum_{\beta =\gamma +\delta }c_{\gamma
,\delta }\left( u^{\prime }\right) ^{\gamma }v^{\delta },
\end{equation*}
we have
\begin{eqnarray*}
\int\limits_{\mathbb{R}}u^{\alpha }\psi (z,u)du &=&\int_{\mathbb{R}
}u^{\alpha }\left\{ \int_{\mathbb{R}}\psi ^{(2)}(u-v)\psi
^{(1)}(z,v)dv\right\} du \\
&=&\int_{\mathbb{R}}\left\{ \int_{\mathbb{R}}\left( u^{\prime }+v\right)
^{\alpha }\psi ^{(2)}(u^{\prime })du\right\} \psi ^{(1)}(z,v)dv \\
&=&\sum_{\alpha =\gamma +\delta }c_{\gamma ,\delta }\int_{\mathbb{R}}\left\{
\int_{\mathbb{R}}\left( u^{\prime }\right) ^{\gamma }\psi ^{(2)}(u^{\prime
})du^{\prime }\right\} v^{\delta }\psi ^{(1)}(z,v)dv \\
&=&\sum_{\alpha =\gamma +\delta }c_{\gamma ,\delta }\int_{\mathbb{R}}\left\{
0\right\} v^{\delta }\psi ^{(1)}(z,v)dv\\
&=&0,
\end{eqnarray*}
for all $\alpha \geq 0$ and each $z\in \mathbb{C}^{n}$. Note
that as a consequence the \emph{full} moments $\int\limits_{\mathbb{H}%
^{n}}z^{\alpha }u^{\beta }\psi (z,u)dzdu$ all vanish, but that in general
the partial moments $\int\limits_{\mathbb{C}^{n}}z^{\alpha }\psi (z,u)dz$ do
\emph{not} vanish.

In order to introduce the flag test function space on the Heisenberg $\mathbb H^n,$ we first introduce
the product test function space on $\mathbb H^n\times \mathbb R$ as follows.
\begin{definition}\label {product test space}
The product test function space on $\mathbb H^n\times \mathbb R$ is the collection of all Schwartz functions
$F$ on $\mathbb H^n\times \mathbb R$ such that
\begin{eqnarray*}
\int_{\mathbb{H}^n}z^{\alpha }u^\beta F( z,u,v)dzdu =0
\end{eqnarray*}
for all $|\alpha|, \beta \geq 0, v\in \mathbb R$ and
\begin{eqnarray*}
\int_{\mathbb{R}}v^{\gamma }F( z,u,v)dv =0
\end{eqnarray*}
for all $\gamma\geq 0, (z,u)\in \mathbb H^n.$
\end{definition}
If $F$ is a product test function on $\mathbb H^n\times \mathbb R,$ we denote $F\in \mathcal S_\infty(\mathbb H^n\times \mathbb R).$
The norm of $F$ is defined by the Schwartz norm of $F$ on $\mathbb H^n\times \mathbb R$ and is denoted by $\|F\|_{\mathcal S_\infty}.$

Applying the projection $\pi$ defined on $\mathbb{H}^{n}$ as introduced in \cite{MRS}:
\begin{equation}
f\left( z,u\right) =\left( \pi F\right) \left( z,u\right) \equiv \int_{
\mathbb{R}}F\left( \left( z,u-v\right) ,v\right) dv,  \label{projection}
\end{equation}
we now define the \emph{flag} test function space $\mathcal{S}_{flag}\left( \mathbb{H}^{n}\right) $ by
\begin{definition}\label {flag test space}
The \emph{flag} test function space $\mathcal{S}_{flag}\left( \mathbb{H}^{n}\right) $ consists of all functions $f(z,u)$
defined on $\mathbb{H}^{n}$ satisfying
\begin{equation}
f\left( z,u\right) =\left( \pi F\right) \left( z,u\right) \equiv \int_{
\mathbb{R}}F\left( \left( z,u-v\right) ,v\right) dv  \label{projection}
\end{equation}
for some $F\in \mathcal S_\infty(\mathbb H^n\times \mathbb R).$ The norm of $f$ is defined by
$$\|f\|_{\mathcal{S}_{flag}}=\inf\{\|F\|_{\mathcal S_\infty}: {\rm for\ all\ } f: f=\pi F\}.$$
\end{definition}
The distribution space, $\mathcal{S}_{flag}^\prime\left(\mathbb{H}^{n}\right),$ is defined by the dual of $\mathcal{S}_{flag}\left( \mathbb{H}^{n}\right).$

Observing that $\psi _{s,t}\left( z,u\right) =\psi _{s}^{(1)}\ast _{2}\psi _{t}^{(2)}\left(
z,u\right)\in \mathcal{S}_{flag}.$ We will show that the Carlder\'on reproducing formula (\ref{CRF st}) also converges in both flag test function space $\mathcal{S}_{flag}\left( \mathbb{H}^{n}\right)$ and distribution space $\mathcal{S}_{flag}^\prime\left( \mathbb{H}^{n}\right)$ as follows.
\begin{theorem} \label{mainCRF}
Let $\psi _{s,t}\left( z,u\right) =\psi _{s}^{(1)}\ast _{2}\psi _{t}^{(2)}\left(
z,u\right).$ Then
\begin{equation}\label{cr}
f(z,u)=\int_{0}^{\infty }\int_{0}^{\infty }{\psi}_{s,t}\ast \psi
_{s,t}\ast f(z,u)\frac{ds}{s}\frac{dt}{t},
\end{equation}
where the integrals converge in both $\mathcal{S}_{flag}( \mathbb{H}^{n})$ and $\mathcal{S}_{flag}^\prime( \mathbb{H}^{n}).$
\end{theorem}
We characterize the flag Lipschitz space $Lip_{flag}^\alpha$ by the following

\begin{theorem}\label{r1} $f\in Lip_{flag}^\alpha$ with $\alpha=(\alpha_1,\alpha_2), \alpha_1, \alpha_2>0$
if and only if $f\in{(\mathcal{S}_{flag})^\prime\left( \mathbb{H}^{n}\right)}$ and
$$\sup_{s,t>0, {(z,u)\in \mathbb H^n}}s^{-\alpha_1}t^{-\alpha_2}|\psi_{s,t}\ast f(z,u)|\leq C<\infty.$$
Moreover,
$$||f||_{Lip_{flag}^\alpha}\sim \sup_{s,t>0,{(z,u)\in \mathbb H^n}}s^{-\alpha_1}t^{-\alpha_2}|\psi_{s,t}\ast f(z,u)|.$$
\end{theorem}

The main result of this paper is

\begin{theorem}\label{r2}
The flag singular integral operator $T$ is bounded on $Lip_{flag}^\alpha$ with $\alpha=(\alpha_1,\alpha_2), \alpha_1,$ $ \alpha_2>0.$
Moreover,
$$||Tf||_{Lip_{flag}^\alpha}\leq C ||f||_{Lip_{flag}^\alpha}.$$
\end{theorem}
As a consequence, the Marcinkiewicz multipliers are bounded on the flag Lischitz space $Lip_{flag}^\alpha.$

\section{Proof of theorem \ref{mainCRF}}

Observe that for $F_1, F_2\in \mathcal S_\infty(\mathbb H^n\times \mathbb R),$ then $\pi(F_1\ast F_2)=(\pi F_1)\ast (\pi F_2).$
Suppose that $f\in {\mathcal S}_{flag}(\mathbb H^n)$ and $f(x, y) =\pi (F)$ with $F\in \mathcal S_\infty(\mathbb H^n\times \mathbb R).$ Set $\Psi_{s,t}(z,u,r)=\psi_s^{(1)}(z,u)\psi_t^{(2)}(r)$ then $\psi_{s,t}=
\pi (\Psi_{s,t}).$ We rewrite the Calder\'on reproducing formula in \eqref{cr} as
$$ (\pi F)(z,u)=\int_{0}^{\infty }\int_{0}^{\infty }\pi({\Psi}_{s,t})\ast ({\pi\Psi}
_{s,t})\ast ({\pi F})(z,u)\frac{ds}{s}\frac{dt}{t}.$$
To show Theorem \ref{mainCRF}, it suffices to prove that for $f\in {\mathcal S}_{flag}(\mathbb H^n),$
$$\bigg\|\pi\big( F -\int_{N^{-1}}^{N }\int_{N^{-1}}^{N }{\Psi}_{s,t}\ast \Psi
_{s,t}\ast F\frac{ds}{s}\frac{dt}{t}\big)\bigg\|_{\mathcal S_{flag}}\rightarrow 0$$
as $N$ tend to $\infty.$

By the definition of the norm of ${\mathcal S}_{flag}(\mathbb H^n),$ one only needs to show that
$$\bigg\|F -\int_{N^{-1}}^{N }\int_{N^{-1}}^{N }{\Psi}_{s,t}\ast \Psi
_{s,t}\ast F(z,u,r)\frac{ds}{s}\frac{dt}{t}\bigg\|_{\mathcal S_\infty}\rightarrow 0$$
as $N$ tend to $\infty.$

To do this, observe that
\begin{eqnarray*}
{\Psi}_{s,t}\ast _{\mathbb{H}^{n}\times \mathbb R}\Psi _{s,t}\ast _{\mathbb{H}
^{n}\times \mathbb R} F\left( (z,u),r\right) &=&\left( {\psi}_{s}^{(1)}\psi
_{t}^{(2)}\right) \ast _{\mathbb{H}^{n}\times \mathbb R}\left( {\psi}_{s}^{(1)}
\psi _{t}^{(2)}\right) \ast _{\mathbb{H}^{n}\times \mathbb R}F\left( (z,u),r\right) \\
&=&\bigg[\left( {\psi}_{s}^{(1)}\ast _{\mathbb{H}^{n}}{\psi}
_{s}^{(1)}\right)  \left(\psi _{t}^{(2)}\ast _{\mathbb{
R}}\psi _{t}^{(2)}\right)\bigg] \ast _{\mathbb{H}^{n}\times \mathbb R}F\left( (z,u),r\right)
\end{eqnarray*}
yields:
\begin{eqnarray*}
&&\int_{0}^{\infty }\int_{0}^{\infty }{\Psi}_{s,t}\ast _{\mathbb{H}
^{n}\times \mathbb R}\Psi _{s,t}\ast _{\mathbb{H}^{n}\times\mathbb R}F((z,u),r)\frac{ds}{s}\frac{dt}{t} \\
&=&\int_{0}^{\infty }{\psi}_{s}^{(1)}\ast _{\mathbb{H}^{n}}{\psi}
_{s}^{(1)}\ast _{\mathbb{H}^{n}}\left\{ \int_{0}^{\infty }\psi
_{t}^{(2)}\ast _{\mathbb{R}}\psi _{t}^{(2)}\ast _{\mathbb R}F\left( (z,u),r\right) \frac{
dt}{t}\right\} \frac{ds}{s} \\
&=&\int_{0}^{\infty }{\psi}_{s}^{(1)}\ast _{\mathbb{H}^{n}}{\psi}
_{s}^{(1)}\ast _{\mathbb{H}^{n}}F\left( (z,u),r\right) \frac{ds}{s}\\&=&F\left(
(z,u),r\right).
\end{eqnarray*}
Write
\begin{align*}
&\bigg|F -\int_{N^{-1}}^{N }\int_{N^{-1}}^{N }{\Psi}_{s,t}\ast \Psi
_{s,t}\ast F(z,u,r)\frac{ds}{s}\frac{dt}{t}\bigg|\\
\le & \int_0^{\infty }\int_0^{N^{-1}}|{\Psi}_{s,t}\ast \Psi
_{s,t}\ast F(z,u,r)|\frac{ds}{s}\frac{dt}{t}+ \int_0^{\infty }\int_N^\infty|{\Psi}_{s,t}\ast \Psi
_{s,t}\ast F(z,u,r)|\frac{ds}{s}\frac{dt}{t}\\&+\int_0^{N^{-1}}\int_0^{\infty }|{\Psi}_{s,t}\ast \Psi
_{s,t}\ast F(z,u,r)|\frac{ds}{s}\frac{dt}{t}+ \int_{N}^{\infty }\int_0^{\infty }|{\Psi}_{s,t}\ast \Psi
_{s,t}\ast F(z,u,r)|\frac{ds}{s}\frac{dt}{t}\\
=& I + II + III + IV.
\end{align*}
We claim that if $N\geq 1$ and for any given an positive integer $L,$ there exists a constant $C_L$ independent of $N$ such that
$$I\leq C_L (1+|(z,u)|^2+ |r|)^{-L}N^{-1}\|F\|_{\mathcal S_{\infty}}.$$
To show the claim, we write
\begin{eqnarray*}
&&{\Psi}_{s,t}\ast _{\mathbb R\times\mathbb{H}^{n}}\Psi _{s,t}\ast _{\mathbb{H}
^{n}\times \mathbb R} F\left( (z,u),r\right)
 \\
 &&=\bigg[\left( {\psi}_{s}^{(1)}\ast _{\mathbb{H}^{n}}{\psi}
_{s}^{(1)}\right)  \left(\psi _{t}^{(2)}\ast _{\mathbb{
R}}\psi _{t}^{(2)}\right)\bigg] \ast _{\mathbb{H}^{n}\times\mathbb R}F\left( (z,u),r\right)\\
&&=\int_{\mathbb R\times\mathbb H^n}\psi_{s}^{1}(z',u')\psi _{t}^{2}(r')F((z,u)\circ(z',u')^{-1},r-r')dz'du'dr',
\end{eqnarray*}
where $\psi_{s}^{1}={\psi}_{s}^{(1)}\ast _{\mathbb{H}^{n}}{\psi}
_{s}^{(1)}$ and $\psi _{t}^{2}=\psi _{t}^{(2)}\ast _{\mathbb{
R}}\psi _{t}^{(2)}.$

Considering the case with $0<t\leq 1$ first and applying the cancellation conditions on $\psi_{s}^{1}$ and $\psi_{t}^{2}$ imply that
\begin{eqnarray}\label{case I}
&&\int_{\mathbb R\times\mathbb H^n}\psi_{s}^{1}(z',u')\psi _{t}^{2}(r')F((z,u)\circ(z',u')^{-1},r-r')dz'du'dr'\nonumber\\
&=&\int_{\mathbb R\times\mathbb H^n}\psi_{s}^{1}(z',u')\psi _{t}^{2}(r')
[F((z,u)\circ(z',u')^{-1},r-r')-F((z,u),r-r')\\
&&-F((z,u)\circ(z',u')^{-1},r)
+F((z,u),r)]dz'du'dr'\nonumber.
\end{eqnarray}
Noting $|(z,u)|^2=|z|^2+|u|$ and applying the smoothness conditions on $F,$ we obtain the following estimates: For any given positive integer $L$ there exists a constant $C$ such that
 \begin{enumerate}
\item[(1)]$|F((z,u)\circ(z',u')^{-1},r-r')-F((z,u),r-r')-F((z,u)\circ(z',u')^{-1},r)+F((z,u),r)|\\\leq C \frac{\displaystyle |(z',u')|}{\displaystyle(1+|(z,u)|)^{L+2n+2}}\frac{\displaystyle |r'|}{\displaystyle(1+|r|)^{L+1}}\|F\|_{\mathcal S_{\infty}},$ \\[5pt]
    for $|(z',u')|\leq \frac{1}{2\gamma}(1+|(z,u)|)$ and $|r'|\leq \frac{1}{2}(1+|r|);$\\

\item[(2)]$|F((z,u)\circ(z',u')^{-1},r-r')-F((z,u),r-r')-F((z,u)\circ(z',u')^{-1},r)+F((z,u),r)|\\\leq C  \frac{\displaystyle |(z',u')|}{\displaystyle(1+|(z,u)|)^{L+2n+2}}\bigg[\frac{\displaystyle1}{\displaystyle(1+|r-r'|)^{L+1}}+\frac{\displaystyle1}{\displaystyle(1+|r|)^{L+1}}\bigg]\|F\|_{\mathcal S_{\infty}}$\\[5pt]
    for $|(z',u')|\leq \frac{1}{2\gamma}(1+|(z,u)|)$ and $|r'|\geq \frac{1}{2}(1+|r|);$\\

\item[(3)]$|F((z,u)\circ(z',u')^{-1},r-r')-F((z,u),r-r')-F((z,u)\circ(z',u')^{-1},r)+F((z,u),r)|\\ \leq C  \bigg(\frac{\displaystyle1}{\displaystyle(1+|(z,u)\circ(z',u')^{-1}|)^{L+2n+2}}+\frac{\displaystyle1}{\displaystyle(1+|(z,u)|)^{L+2n+2}}\bigg)\frac{\displaystyle |r'|}{\displaystyle(1+|r|)^{L+1}}\|F\|_{\mathcal S_{\infty}}$\\[5pt]
    for $|(z',u')|\geq \frac{1}{2\gamma}(1+|(z,u)|)$ and $|r'|\leq \frac{1}{2}(1+|r|);$\\

\item[(4)]$|F((z,u)\circ(z',u')^{-1},r-r')-F((z,u),r-r')-F((z,u)\circ(z',u')^{-1},r)+F((z,u),r)|\\\leq C \bigg (\frac{\displaystyle1}{\displaystyle(1+|(z,u)\circ(z',u')^{-1}|)^{L+2n+2}}+\frac{\displaystyle1}{\displaystyle(1+|(z,u)|)^{L+2n+2}}\bigg)\\  \,
 \hskip1cm\times\bigg(\frac{\displaystyle1}{\displaystyle(1+|r-r'|)^{L+1}}+\frac{\displaystyle1}{\displaystyle(1+|r|)^{L+1}}\bigg)\|F\|_{\mathcal S_{\infty}}$\\[5pt]
    for $|(z',u')|\geq \frac{1}{2\gamma}(1+|(z,u)|)$ and $|r'|\geq \frac{1}{2}(1+|r|).$
\end{enumerate}
Plugging the above estimates together with the size conditions
\begin{eqnarray*}
&|\psi_{s}^{1}(z',u')|\leq C s^{-2n-2}\frac{\displaystyle 1}{\displaystyle \bigg(1+\big|\frac{z'}{s}\big|^2+\big|\frac{u'}{s^2}\big|\bigg)^{n+2+L}},\\
&|\psi _{t}^{2}(r')|\leq C t^{-1}\frac{\displaystyle 1}{\displaystyle \Big(1+\big|\frac{r'}{t}\big|\Big)^{2+L}}
\end{eqnarray*}

into \eqref{case I} implies that whenever $0<t\leq 1,$
\begin{eqnarray*}
&&\bigg|\int_{\mathbb R\times\mathbb H^n}\psi_{s}^{1}(z',u')\psi _{t}^{2}(r')F((z,u)\circ(z',u')^{-1},r-r')dz'du'dr'\bigg|\\
&\leq & C{s}{t}\frac{1}{(1+|(z,u)|^2+|r|)^{L}}\|F\|_{\mathcal S_{\infty}}.
\end{eqnarray*}
We now consider the case with $1\leq t<\infty.$
For this case, applying the cancellation conditions on $\psi_{s}^{1}$ and $F((z,u),r)$ with respect to the variable $r,$ we write
\begin{eqnarray}\label{case II}
&&\int_{\mathbb R\times\mathbb H^n}\psi_{s}^{1}(z',u')\psi _{t}^{2}(r')F((z,u)\circ(z',u')^{-1},r-r')dz'du'dr'\nonumber\\
&=&\int_{\mathbb R\times\mathbb H^n}\psi_{s}^{1}(z',u')\Big[\psi _{t}^{2}(r') -\sum_{j\leq L}c_j(\psi _{t}^{2})^{(j)}(r)(r'-r)^j\Big]\\
&&\times
\Big[F((z,u)\circ(z',u')^{-1},r-r')-F((z,u),r-r')\Big]dz'du'dr'\nonumber,
\end{eqnarray}
where $\sum_{j\leq L}c_j(\psi _{t}^{2})^{(j)}(r)(r'-r)^j$ is the Taylor series of $\psi _{t}^{2}$ at $r$ with the degree $L.$
Inserting estimates
 \begin{enumerate}
\item[(1)] $|\psi _{t}^{2}(r') -\sum_{j\leq L}c_j(\psi _{t}^{2})^{(j)}(r)(r'-r)^j|\leq C\bigg(\frac{\displaystyle |r-r'|}{\displaystyle t+|r|}\bigg)^{L}\frac{\displaystyle t^L}{\displaystyle(t+|r|)^{1+L}},$\\
for $|r-r'|\leq \frac{1}{2}(t+|r|);$\\

\item[(2)] $|\psi _{t}^{2}(r') -\sum_{j\leq L}c_j(\psi _{t}^{2})^{(j)}(r)(r'-r)^j|\leq C\bigg[\frac{\displaystyle t^L}{\displaystyle(t+|r'|)^{1+L}}+\sum_{j\leq L}c_j\frac{\displaystyle |r'-r|^j}{\displaystyle t^j}\frac{\displaystyle t^L}{\displaystyle(t+|r|)^{1+L}}\bigg]$\\
for $|r-r'|\geq \frac{1}{2}(t+|r|)$
\end{enumerate}
together with
\begin{enumerate}
\item[(3)] $|F((z,u)\circ(z',u')^{-1},r-r')-F((z,u),r-r')|\\ \leq C
\frac{\displaystyle |(z',u')|}{\displaystyle (1+|(z,u)|)^{L+2n+2}} \frac{\displaystyle1}{\displaystyle(1+|r-r'|)^{L+2}}\|F\|_{\mathcal S_{\infty}},$\\
for $|(z',u')|\leq \frac{1}{2\gamma}(1+|(z,u)|)$;\\

\item[(4)] $|F((z,u)\circ(z',u')^{-1},r-r')-F((z,u),r-r')|\\ \leq C
\bigg(\frac{\displaystyle1}{\displaystyle(1+|(z,u)\circ(z',u')^{-1}|)^{L+2n+2}}+\frac{\displaystyle1}{\displaystyle(1+|(z,u)|)^{L+2n+2}}\bigg)\frac{\displaystyle1}{\displaystyle(1+|r-r'|)^{L+1}}\|F\|_{\mathcal S_{\infty}}$\\
for $|(z',u')|\geq \frac{1}{2\gamma}(1+|(z,u)|)$
\end{enumerate}
into \eqref{case II} gives that if $t\geq 1,$
\begin{eqnarray*}
&&\bigg|\int_{\mathbb R\times\mathbb H^n}\psi_{s}^{1}(z',u')\psi _{t}^{2}(r')F((z,u)\circ(z',u')^{-1},r-r')dz'du'dr'\bigg|\\
&\leq & C\frac{s}{t}\frac{1}{(1+|(z,u)|^2+|r|)^{L}}\|F\|_{\mathcal S_{\infty}}.
\end{eqnarray*}
These estimates imply the proof of the claim.
Note that $(\partial_z^\alpha\partial_u^\beta\partial_r^\gamma F)$ satisfy similar regularity and cancellation conditions as $F$ does, and
\begin{eqnarray*}
&\partial_z^\alpha\partial_u^\beta\partial_r^\gamma\check{\Psi}_{s,t}\ast _{\mathbb{H}^{n}\times \mathbb R}\Psi _{s,t}\ast _{\mathbb{H}
^{n}\times \mathbb R} F\left( (z,u),r\right)\\
&=\check{\Psi}_{s,t}\ast _{\mathbb{H}^{n}\times \mathbb R}\Psi _{s,t}\ast _{\mathbb{H}
^{n}\times \mathbb R} (\partial_z^\alpha\partial_u^\beta\partial_r^\gamma F)\left( (z,u),r\right).
\end{eqnarray*}
These facts yield that
$$\|I\|_{\mathcal S_{\infty}}\leq CN^{-1}$$
and hence $\|I\|_{\mathcal S_{\infty}}\rightarrow 0$ as $N$ tends to $\infty.$ The proofs for $II, III$ and $IV$ are similar and thus the proof of Theorem \ref{mainCRF} is concluded.

\section{Proof of theorem \ref{r1}}

We first show that if $f\in Lip_{flag}^\alpha$ with $\alpha=(\alpha_1,\alpha_2), 0<\alpha_1, \alpha_2<1,$ then $f\in {\mathcal S}_{flag}^\prime.$ To see this, for each $g(z,u)\in \mathcal{S}_{flag}$ where $g=\pi G$ with $G\in \mathcal S_\infty$ and $\|G\|_{\mathcal S_\infty}=1,$ we observe that
\begin{eqnarray*}
\langle g,f\rangle &=&\int_{\mathbb H^n}\int_{\mathbb R}G((z,u-v),v)f(z,u)dvdzdu\\
&=&\int_{\mathbb H^n}\int_{\mathbb R}G((z,u),v)f(z,u+v)dvdzdu.
\end{eqnarray*}
Applying the cancellation conditions on $G((z,u),v),$ we have
\begin{eqnarray*}
\langle g,f\rangle =\int_{\mathbb H^n}\int_{\mathbb R}G((z,u),v)[f(z,u+v)-f(0,v)-f(z,u)+f(0,0)]dvdzdu.
\end{eqnarray*}
The size conditions on $G$ and the fact that $f\in Lip_{flag}^\alpha$ imply that
\begin{eqnarray}\label{31}
|\langle g,f\rangle |\leq C \|f\|_{Lip^\alpha_{flag}},
\end{eqnarray}
which shows that $f\in (\mathcal S_{flag})'.$

Observe that

$$\psi_{s,t}\ast f(z,u)= \int_{\mathbb H^n}\int_{\mathbb R}\psi_s^{(1)}(z',u')\psi_t^{(2)}(v)f((z,u)\circ (z',u'+v)^{-1})dvdz'du'.$$

Applying the cancellations conditions on $\psi_s^{(1)}$ and $\psi_t^{(2)},$ we have
\begin{eqnarray*}
\psi_{s,t}\ast f(z,u)
&=&\int_{\mathbb H^n}\int_{\mathbb R}\psi_s^{(1)}(z',u')\psi_t^{(2)}(v)[f((z,u)\circ(z',u'+v)^{-1})\\
&&-f((z,u)\circ(z',u')^{-1})-f((z,u)\circ(0,v)^{-1})+f(z,u)]dvdz'du'\\
&=&\int_{\mathbb H^n}\int_{\mathbb
R}\psi_s^{(1)}(z',u')\psi_t^{(2)}(v)\Delta^2_{v}\Delta^{1}_{(z',u')}(f)(z,u)dvdz'du'.
\end{eqnarray*}
The size conditions on $\psi_s^{(1)}$ and $\psi_t^{(2)}$ and the fact that $f\in Lip_{flag}^\alpha$ yield that
\begin{eqnarray}\label{32}
|\psi_{s,t}\ast f(z,u)|\leq C s^{\alpha_1}t^{\alpha_2}\|f\|_{Lip^\alpha_{flag}},
\end{eqnarray}
which implies that
$s^{-\alpha_1}t^{-\alpha_2}|\psi_{s,t}\ast f(z,u)|\leq C ||f||_{Lip_{flag}^\alpha}$ for any $s,t>0,$ and $(z,u) \in \mathbb H^n$.

When $\alpha=(\alpha_1,\alpha_2)$ with
$\alpha_1=1, 0<\alpha_2<1,$ note first that $\psi_s^{(1)}$ is a radial function and then applying the cancellation conditions on $\psi_t^{(2)},$ we have
\begin{eqnarray*}
&&\psi_{s,t}\ast f(z,u)\\&=&\frac{1}{2}\int_{\mathbb H^n}\int_{\mathbb R}\psi_s^{(1)}(z',u')\psi_t^{(2)}(v)[f((z,u)\circ (z',u'-v))\\
&&+f((z,u)\circ (z',u'+v)^{-1})]dvdz'du'\\
&=&\frac{1}{2}\int_{\mathbb H^n}\int_{\mathbb R}\psi_s^{(1)}(z',u')\psi_t^{(2)}(v)[f((z,u)\circ (z',u'-v))\\
&&+f((z,u)\circ (z',u'+v)^{-1})-2f((z,u)\circ (0,v)^{-1})]dvdz'du'
\\&=&\frac{1}{2}\int_{\mathbb H^n}\int_{\mathbb R}\psi_s^{(1)}(z',u')\psi_t^{(2)}(v)\{[f((z,u)\circ (z',u'-v))+f((z,u)\circ (z',u'+v)^{-1})\\
&&-2f((z,u)\circ (0,v)^{-1})]\nonumber -[f((z,u)\circ
(z',u'))+f((z,u)\circ (z',u')^{-1})-2f(z,u)]\}dvdz'du' \\
&=&\frac{1}{2}\int_{\mathbb H^n}\int_{\mathbb R}\psi_s^{(1)}(z',u')\psi_t^{(2)}(v)\Delta^2_{v}\Delta^{1,Z}_{(z',u')}(f)(z,u)dvdz'du'
.
\end{eqnarray*}
Applying the same proof gives the desired estimate for this case.
All other cases $\alpha=(\alpha_1,\alpha_2)$ where
$0<\alpha_1<1, \alpha_2=1$
or $\alpha_1=\alpha_2=1$ can be handled similarly.
For the case where $\alpha=(\alpha_1,\alpha_2)$ with
$\alpha_i>1, \alpha_i=m_i+r_i, 0<r_i\leq 1, i=1,2,$ set ${\partial^{m_1+m_2}\widetilde\psi}(z,u)=\psi$ with  $s^{-m_1}t^{-m_2}{\widetilde \psi}_{s,t}$ satisfy the similar smoothness, size and cancellation conditions as $\psi_{s,t}.$ Thus, repeating the same proof gives
$$|\psi_{s,t}\ast f|=|s^{m_1}t^{m_2}(s^{-m_1}t^{-m_2}{\widetilde \psi}_{s,t})\ast\partial^{m_1+m_2}f|\leq Cs^{m_1}t^{m_2}s^{r_1}t^{r_2}\|f\|_{Lip^\alpha_{flag}}=Cs^{\alpha_1}t^{\alpha_2}\|f\|_{Lip^\alpha_{flag}}.$$
Therefore, this case can be also handled similarly.

We now prove the converse implication of Theorem \ref{r1}. Suppose that $f\in \mathcal S_{flag}^\prime(\mathbb H^n)$ and $|\psi_{s,t}\ast f(z,u)|\leq Cs^{\alpha_1}t^{\alpha_2}$ with $\alpha_1, \alpha_2>0.$ We first show that $f$ is a continuous function on $\mathbb H^n.$ To do this, by Theorem \ref{mainCRF}, $$f(z,u)=\int^\infty_0\int^\infty_0 \psi_{s,t}\ast\psi_{s,t}\ast f(z,u)\frac{ds}{s}\frac{dt}{t},$$ where the integral converges in $\mathcal S_{flag}^\prime(\mathbb H^n)$. We split the above integral by the following four cases:
$$(i)\ s,t<1;\ (ii)\ s\geq 1, t<1;\ (iii)\ s<1, t\geq 1; {\rm\ and\ }(iv)\ s,t\geq 1.$$
Then
we write $f=f_1+f_2+f_3+f_4$ in $\mathcal S_{flag}^\prime(\mathbb H^n)$ according to the above four cases of $s$ and $t$. As for $f_1$, observe that the fact $|\psi_{s,t}\ast\psi_{s,t}\ast f(z,u)|\leq Cs^{\alpha_1}t^{\alpha_2}$ implies the integral $$\int^1_0\int^1_0 \psi_{s,t}\ast\psi_{s,t}\ast f(z,u)\frac{ds}{s}\frac{dt}{t}$$ converges uniformly and hence, $f_1$ is a continuous function on $\mathbb H^n.$

To see that $f_2$ is continuous, for $g\in \mathcal S_{flag},$
$$\langle f_2, g\rangle =\int_0^1\int_1^\infty\langle \psi_{s,t}\ast\psi_{s,t}\ast f(z,u), g(z,u)\rangle \frac{ds}{s}\frac{dt}{t}.$$
Note that by the cancellation condition on $g$ and for any fixed positive integer $N,$ we have
\begin{eqnarray*}
&&\langle \psi_{s,t}\ast\psi_{s,t}\ast f(z,u), g(z,u)\rangle \\
&=&\bigg\langle \int_{\mathbb H^n}\int_{\mathbb R}\psi_s^{(1)}((z,u)\circ (z',u'+v)^{-1})\psi_t^{(2)}(v)\psi_{s,t}\ast f(z',u')dvdz'du', g(z,u)\bigg\rangle \\
&=&\bigg\langle \int_{\mathbb H^n}\int_{\mathbb R}\Big[\psi_s^{(1)}((z,u)\circ (z',u'+v)^{-1})-\sum_{|\alpha|+2\beta\leq 2N}c_{\alpha,\beta}\partial^\alpha_z\partial^\beta_u \psi_s^{(1)}((z',u'+v)^{-1})(z,u)^{\alpha, \beta}\Big]\\
&&\times\psi_t^{(2)}(v)\psi_{s,t}\ast f(z',u')dvdz'du', g(z,u)\bigg\rangle ,
\end{eqnarray*}
which implies that
\begin{eqnarray*}
f_2(z,u)&=&\int_0^1\int_1^\infty\int_{\mathbb H^n}\int_{\mathbb R}\Big[\psi_s^{(1)}((z,u)\circ (z',u'+v)^{-1})\\
&&-\sum_{|\alpha|+2\beta\leq 2N}c_{\alpha,\beta}\partial^\alpha_z\partial^\beta_u \psi_s^{(1)}((z',u'+v)^{-1})(z,u)^{\alpha, \beta}\Big]
\psi_t^{(2)}(v)\psi_{s,t}\ast f(z',u')dvdz'du'\frac{ds}{s}\frac{dt}{t}
\end{eqnarray*}
in the sense of $\mathcal S_{flag}^\prime(\mathbb H^n).$

Observe that
\begin{eqnarray*}
&&\bigg|\int_{\mathbb H^n}\int_{\mathbb R}\Big[\psi_s^{(1)}((z,u)\circ (z',u'+v)^{-1})\\
&&-\sum_{|\alpha|+2\beta\leq 2N}c_{\alpha,\beta}\partial^\alpha_z\partial^\beta_u \psi_s^{(1)}((z',u'+v)^{-1})(z,u)^{\alpha, \beta}\Big]\psi_t^{(2)}(v)\psi_{s,t}\ast f(z',u')dvdz'du'\bigg|\\
&&\leq Cs^{-2N+\alpha_1}t^{\alpha_2}|(z,u)|^{2N}.
\end{eqnarray*}
The above estimate implies that for any given $R>0,$ the integral for $f_2$ is converges uniformly for $|(z,u)|^{2N}\leq R$ and thus, $f_2$ is a continuous function on any compact subset in $\mathbb H^n.$

Note that $(z,u)\circ (z', u'+v)^{-1}=(z, -u')\circ (z', -u+v)^{-1}.$ We write
\begin{eqnarray*}
&&\langle \psi_{s,t}\ast\psi_{s,t}\ast f(z,u), g(z,u)\rangle \\
&=&\bigg\langle \int_{\mathbb H^n}\int_{\mathbb R}\psi_s^{(1)}((z,-u')\circ (z',-u+v)^{-1})\psi_t^{(2)}(v)\psi_{s,t}\ast f(z',u')dvdz'du', g(z,u)\bigg\rangle \\
&=&\bigg\langle \int_{\mathbb H^n}\int_{\mathbb R}\psi_s^{(1)}((z,-u')\circ (z',v)^{-1})\psi_t^{(2)}(u+v)\psi_{s,t}\ast f(z',u')dvdz'du', g(z,u)\bigg\rangle \\
&=& \bigg\langle \int_{\mathbb H^n}\int_{\mathbb R}\psi_s^{(1)}((z,-u')\circ (z',v)^{-1})\Big[\psi_t^{(2)}(u+v)-\sum_{\gamma \leq 2N}c_{\gamma }\frac{d^\gamma}{dv^\gamma}\psi_t^{(2)}(v)v^\gamma\Big]\\
&&\times\psi_{s,t}\ast f(z',u')dvdz'du', g(z,u)\bigg\rangle .
\end{eqnarray*}
Similarly,
\begin{eqnarray*}
&& \bigg|\int_{\mathbb H^n}\int_{\mathbb R}\psi_s^{(1)}((z,-u')\circ (z',v)^{-1})\Big[\psi_t^{(2)}(u+v)-\sum_{\gamma \leq 2N}c_{\gamma }\frac{d^\gamma}{dv^\gamma}\psi_t^{(2)}(v)v^\gamma\Big]\psi_{s,t}\ast f(z',u')dvdz'du'\bigg|\\
&&\leq C s^{\alpha_1}t^{-2N+\alpha_2}|u|^{2N}
\end{eqnarray*}
and hence $f_3$ is a continuous function on any compact subset in $\mathbb H^n.$

Taking the geometric means of these two estimates shows that $f_4$ is a continuous function.

Now we show that $ f\in Lip_{flag}^\alpha.$ First, if $\alpha=(\alpha_1, \alpha_2)$ with $0<\alpha_1, \alpha_2<1,$ then
\begin{eqnarray*}&&\Big|\Delta^2_{w}\Delta^1_{(u,v)}(f)(z,r)\Big|\\
&=&\bigg|f((z,r)\circ(u,v+w)^{-1})-f((z,r)\circ(u,v)^{-1})-f((z,r)\circ(0,w)^{-1})+f(z,r)\bigg|
\\&=&\Big|\int_0^\infty\int_0^\infty\int_{\mathbb H^n}\Big[\psi_{s,t}((z,r)\circ(u,v+w)^{-1}\circ (z',r')^{-1})-\psi_{s,t}((z,r)\circ(u,v)^{-1}\circ(z',r'))\\
&&-\psi_{s,t}((z,r)\circ(0,w)^{-1}\circ(z',r')^{-1})+\psi_{s,t}((z,r)\circ(z',r'))\Big]
\psi_{s,t}\ast f(z',r')dz'dr'\frac{ds}{s}\frac{dt}{t}\Big|.
\end{eqnarray*}
Observe that
\begin{eqnarray*}
&&\psi_{s,t}((z,r)\circ(u,v+w)^{-1}\circ (z',r')^{-1})-\psi_{s,t}((z,r)\circ(u,v)^{-1}\circ(z',r'))\\
&&-\psi_{s,t}((z,r)\circ(0,w)^{-1}\circ(z',r')^{-1})+\psi_{s,t}((z,r)\circ(z',r'))\\
&=&\int_{\mathbb R }\Big[\psi_s^{(1)}((z,r)\circ(u,v)^{-1}\circ(z',v'+r')^{-1})\\
&&-\psi_s^{(1)}((z,r)\circ(z',v'+r')^{-1})\Big]\Big[\psi_t^{(2)}(v'-w)-\psi_t^{(2)}(v')\Big]dv'.
\end{eqnarray*}
We have
\begin{eqnarray*}
A&=&\int_{\mathbb H^n}\Big[\psi_{s,t}((z,r)\circ(u,v+w)^{-1}\circ (z',r')^{-1})-\psi_{s,t}((z,r)\circ(u,v)^{-1}\circ(z',r'))\\
&&-\psi_{s,t}((z,r)\circ(0,w)^{-1}\circ(z',r')^{-1})+\psi_{s,t}((z,r)\circ(z',r'))\Big]
\psi_{s,t}\ast f(z',r')dz'dr'\\
&=&\int_{\mathbb H^n}\int_{\mathbb R }\Big[\psi_s^{(1)}((z,r)\circ(u,v)^{-1}\circ(z',v'+r')^{-1})-\psi_s^{(1)}((z,r)\circ(z',v'+r')^{-1})\Big]\\
&&\times\Big[\psi_t^{(2)}(v'-w)-\psi_t^{(2)}(v')\Big] \psi_{s,t}\ast
f(z',r')dv'dz'dr'.
\end{eqnarray*}
We now choose $t_0$ and $s_0$ such that $t_0\leq |(u,v)|<2t_0$ and
$s_0\leq |w|<2s_0$ and split
\begin{eqnarray*}
\int_0^\infty\int_0^\infty A&=&\int_0^{s_0}\int_0^{t_0} A + \int_{s_0}^{\infty}\int_0^{t_0} A + \int_0^{s_0}\int_{t_0}^{\infty} A +\int_{s_0}^{\infty}\int_{t_0}^{\infty} A\\
&=:&A_1 +A_2+A_3+A_4.
\end{eqnarray*}

To deal with $A_1,$ applying the size conditions on both $\psi^{(1)}_{s}$ and $\psi^{(2)}_{t}$ yields
\begin{eqnarray*}
|A| \leq C|\psi_{s,t}\ast f(z',r')|\leq Cs^{\alpha_1}t^{\alpha_2}, \
\text{for all } s,t>0
\end{eqnarray*}
and hence $ A_1$ is dominated by
$$C\int_0^{s_0}\int_0^{t_0}s^{\alpha_1}t^{\alpha_2}\frac{ds}{s}\frac{dt}{t}\leq Cs_0^{\alpha_1}t_0^{\alpha_2}\leq C|(u,v)|^{\alpha_1}|w|^{\alpha_2}.$$

To estimate $A_2$, applying the smoothness condition on $\psi^{(1)}_{s}$ and the size condition on $\psi^{(2)}_{t}$ implies
\begin{eqnarray*}
|A|\lesssim\frac{|(u,v)|}{s}|\psi_{s,t}\ast f(z',r')|\lesssim
s^{\alpha_1-1}t^{\alpha_2}|(u,v)|.
\end{eqnarray*}
This implies that $A_2$ is bounded by
$$C\int_{s_0}^{\infty}\int_0^{t_0}s^{\alpha_1-1}t^{\alpha_2}\frac{ds}{s}\frac{dt}{t}|(u,v)|\leq Cs_0^{\alpha_1-1}t_0^{\alpha_2}|(u,v)|\leq C|(u,v)|^{\alpha_1}|w|^{\alpha_2}.$$

The estimate for $A_3$ is similar to the estimate for $A_2$. Finally, to handle with $A_4$, applying the smoothness conditions on  both $\psi^{(1)}_{s}$ and $\psi^{(2)}_{t}$ we obtain that
$$|A|\lesssim \frac{|(u,v)|}{s} \cdot \frac{|w|}{t} |\psi_{s,t}\ast f(z',r')|\lesssim s^{\alpha_1-1}t^{\alpha_2-1}|(u,v)||w|.$$
Hence this implies that $A_4$ is dominated by
$$C\int_{s_0}^{\infty}\int_{t_0}^{\infty}s^{\alpha_1-1}t^{\alpha_2-1}\frac{ds}{s}\frac{dt}{t}|(u,v)||w|\leq Cs_0^{\alpha_1-1}t_0^{\alpha_2-1}|(u,v)||w|\leq C|(u,v)|^{\alpha_1}|w|^{\alpha_2}.$$

When $\alpha=(\alpha_1, \alpha_2)$ with $\alpha_1=\alpha_2=1$ observe that
\begin{eqnarray*}&&\Delta^{2,Z}_{w}\Delta^{1,Z}_{(u,v)}(f)(z,r)\\
&=&\Big[f((z,r)\circ(u,v+w))+f((z,r)\circ(u,v-w)^{-1})-2f((z,r)\circ(0,w))\Big]\nonumber\\
&&+\Big[f((z,r)\circ(u,v-w))+f((z,r)\circ(u,v+w)^{-1})-2f((z,r)\circ(0,w)^{-1})\Big]\nonumber\\
&&-2\Big[f((z,r)\circ(u,v))+f((z,r)\circ(u,v)^{-1})-2f(z,r))\Big]
\\&=&\int_{0}^{\infty}\int_{0}^{\infty}\int_{\mathbb H^n}\int_{\mathbb R}\Big[\psi_s^{(1)}((z,r)\circ(u,v)^{-1}\circ(z',v'+r')^{-1})+ \psi_s^{(1)}((z,r)\circ(u,v)\circ(z',v'+r')^{-1}) \\
&&-2\psi_s^{(1)}((z,r)\circ(z',v'+r')^{-1})\Big]
\Big[\psi_t^{(2)}(v'-w)+\psi_t^{(2)}(v'+w) -2\psi_t^{(2)}(v')\Big]
\psi_{j,k}\ast f(z',r')dv'dz'dr'\frac{ds}{s}\frac{dt}{t}.\end{eqnarray*}
Repeating a similar calculation gives the desired result for this case. The other two cases where $\alpha_1=1, 0<\alpha_2<1$ and $0<\alpha_1<1, \alpha_2=1$ can be handled similarly. Lastly, when $\alpha_1=m_1+r_1, \alpha_2=m_2+r_2,$ note that
\begin{eqnarray*}
&&\partial^{m_1+m_2}f((z,r)\circ(u,v+w)^{-1})-\partial^{m_1+m_2}f((z,r)\circ(u,v)^{-1})\\
&&-\partial^{m_1+m_2}f((z,r)\circ(0,w)^{-1})+\partial^{m_1+m_2}f(z,r)
\\&=&\int_{\mathbb H^n}\int_{\mathbb R }\Big[\partial^{m_1}\psi_s^{(1)}((z,r)\circ(u,v)^{-1}\circ(z',v'+r')^{-1})-\partial^{m_1}\psi_s^{(1)}((z,r)\circ(z',v'+r')^{-1})\Big]\\ &&\times
\Big[\partial^{m_2}\psi_t^{(2)}(v'-w)-\partial^{m_2}\psi_t^{(2)}(v')\Big]\psi_{s,t}\ast f(z',r')dv'dz'dr'.
\end{eqnarray*}
Again observe that the properties of $\partial^{m_1}\psi^{(1)}_{s}$ and $\partial^{m_2}\psi^{(2)}_{t}$ are similar to $s^{-m_1}\psi^{(1)}_{s}$ and $t^{-m_2}\psi^{(2)}_{t},$ respectively, and hence the estimate for this case is the same as the proof for the case where $\alpha=(\alpha_1, \alpha_2)$ with $0<\alpha_1, \alpha_2\leq 1.$ We leave the details to the reader. The proof of Theorem \ref{r1} is concluded.

\section{Proof of Theorem \ref{r2}}\label{sec:onbwavelet}

We first show that if $f\in Lip_{flag}^\alpha$ with $\alpha=(\alpha_1, \alpha_2), \alpha_1, \alpha_2>0,$ then there exists a sequence $\{f_n\}$ such that $f_n\in L^2\cap Lip_{flag}^\alpha$ and $f_n$ converges to $f$ in the distribution sense. Moreover,
$||f_n||_{Lip_{flag}^\alpha}\leq C||f||_{Lip_{flag}^\alpha},$ where the constant $C$ is independent of $f_n$ and $f.$ To do this, note that, by Theorem \ref{mainCRF}, for each $f\in Lip_{flag}^\alpha,$
$$f(z,u)=\int_0^\infty\int_0^\infty \psi_{s,t}\ast\psi_{s,t}\ast f(z,u)\frac{ds}{s}\frac{dt}{t},$$
where the integral converges in $\mathcal S_{flag}^\prime(\mathbb H^n)$.
Set
$$f_n=\int_{n^{-1}}^n \int_{n^{-1}}^n \psi_{s,t}\ast\psi_{s,t}\ast f(z,u)\frac{ds}{s}\frac{dt}{t}.$$
Obviously, $f_n\in L^2$ and converges to $f$ in the distribution sense. To see that $f_n\in Lip_{com}^\alpha,$ by Theorem \ref{r1}, $$||f_n||_{Lip_{flag}^\alpha}\leq C\sup\limits_{s,t>0,(z,u)\in \mathbb H^n}s^{-\alpha_1}t^{-\alpha_2}|\psi_{s,t}\ast f_n(z,u)|.$$
Note that
$$\psi_{s,t}\ast f_n(x)=\int_{n^{-1}}^n \int_{n^{-1}}^n\psi_{s,t}\ast\psi_{s',t'}\ast\psi_{s',t'}\ast f(z,u)\frac{ds'}{s'}\frac{dt'}{t'}.$$

By an estimate given in \cite{HLS}, that there exists a constant $C=C(M)$ depending only on $M$ such that
if $\left( {s}\vee {s^{\prime }}\right) ^{2}\leq {t}\vee {t^{\prime }},$
then
\begin{eqnarray*}
&&\left\vert \psi _{s,t}\ast \psi _{s^{\prime },t^{\prime
}}(z,u)\right\vert\\ &\leq &C\left( {\frac{{s}}{{s^{\prime }}}}\wedge {\frac{{
s^{\prime }}}{{s}}}\right) ^{2M}\left( {\frac{{t}}{{t^{\prime }}}}\wedge {
\frac{{t^{\prime }}}{{t}}}\right) ^{M}  {\frac{\left( {{s}\vee {s^{\prime }}}\right) ^{2M}}{\left( {{s}\vee
{s^{\prime }}+}\left\vert {z}\right\vert \right) {^{2n+2M}}}}{\frac{\left( {{
t}\vee {t^{\prime }}}\right) ^{M}}{\left( {{t}\vee {t^{\prime }}+}\left\vert
{u}\right\vert \right) {^{1+M}}}},
\end{eqnarray*}
and if $\left( {s}\vee {s^{\prime }}\right) ^{2}\geq {t}\vee {t^{\prime }},$
then
\begin{eqnarray*}
&&\left\vert \psi _{s,t}\psi _{s^{\prime },t^{\prime }}(z,u)\right\vert\\
&\leq &C\left( {\frac{{s}}{{s^{\prime }}}}\wedge {\frac{{s^{\prime }}}{{s}}}
\right) ^{M}\left( {\frac{{t}}{{t^{\prime }}}}\wedge {\frac{{t^{\prime }}}{{t
}}}\right) ^{M} {\frac{\left( {{s}\vee {s^{\prime }}}\right) ^{M}}{\left( {{s}\vee {
s^{\prime }}+}\left\vert {z}\right\vert \right) {^{2n+M}}}}{\frac{\left( {{s}
\vee {s^{\prime }}}\right) ^{M}}{\left( {{s}\vee {s^{\prime }}+}\sqrt{
\left\vert {u}\right\vert }\right) {^{2+2M}}}}.
\end{eqnarray*}
These estimates imply that if $M>\alpha_1\vee \alpha_2,$
\begin{eqnarray*}
&&s^{-\alpha_1}t^{-\alpha_2}|\psi_{s,t}\ast f_n(z,u)|\\
&\leq& C s^{-\alpha_1}t^{-\alpha_2}\int_0^\infty\int_0^\infty \left( {\frac{{s}}{{s^{\prime }}}}\wedge {\frac{{s^{\prime }}}{{s}}}
\right) ^{M}\left( {\frac{{t}}{{t^{\prime }}}}\wedge {\frac{{t^{\prime }}}{{t
}}}\right) ^{M}(s')^{\alpha_1}(t')^{\alpha_2}\frac{ds'}{s'}\frac{dt'}{t'}\|f\|_{Lip_{flag}}\\
&\leq& C \|f\|_{Lip_{flag}},
\end{eqnarray*}
which implies that $||f_n||_{Lip_{flag}^\alpha}\leq C||f||_{Lip_{flag}^\alpha}.$

Now we claim that if $f\in L^2$ and $T=K\ast f$ is a flag singular integral operator on $\mathbb H^n$ with a flag kernel $K$ as given in \cite{MRS}, then
\begin{eqnarray}\label{3.1}||T(f)||_{Lip_{flag}^\alpha}&\lesssim &||f||_{Lip_{flag}^\alpha}.\end{eqnarray}
Indeed, by  Theorem \ref{r1},
$$||T(f)||_{Lip_{flag}^\alpha}\lesssim  \sup\limits_{s,t>0,(z,u)\in \mathbb H^n}s^{-\alpha_1}t^{-\alpha_2}|\psi_{s,t}\ast T(f)(z,u)|.$$
Observe that, by a result in \cite{MRS}, $T$ is bounded on $L^2(\mathbb H^n)$ and hance
$$\psi_{s,t}\ast T(f)(z,u)=\int_0^\infty\int_0^\infty\psi_{s,t}\ast K\ast\psi_{s',t'}\ast\psi_{s',t'}\ast f(z,u)\frac{ds'}{s'}\frac{dt'}{t'}.$$
Again, By a result in \cite{MRS}, $K(z,u)=\int_{\mathbb R}K^\sharp (z, u-v, v)dv,$ where $K^\sharp(z, u, v), (z,u)\in \mathbb H^n, v\in \mathbb R,$ is a product singular integral kernel on $\mathbb H^n\times \mathbb R.$ Note that
$$\psi_{s,t}\ast K\ast\psi_{s',t'}(z,u)=\int_{\mathbb R}\Psi_{s,t}\ast_{\mathbb H^n\times \mathbb R}K^\sharp\ast_{\mathbb H^n\times \mathbb R}\Psi_{s',t'}(z, u-v, v)dv,$$
where $\pi \Psi_{s,t}=\psi_{s,t}.$

Applying the classical almost orthogonal estimates with $\Psi_{s,t}, K^\sharp$ and $\Psi_{s',t'}$ on $\mathbb H^n\times \mathbb R,$ we have that for any positive integer $M,$
\begin{eqnarray*}
&&|\Psi _{s,t}\ast _{\mathbb{H}^{n}\times \mathbb{R
}}K^\sharp \ast_{\mathbb{H}^{n}\times \mathbb{R
}} \Psi _{s^{\prime },t^{\prime }}\left( ( z,u) ,v\right)|
 \\
&\lesssim &C\left( {\frac{{s}}{{s^{\prime }}}}\wedge {\frac{{s^{\prime }}}{{s
}}}\right) ^{2M}\left( {\frac{{t}}{{t^{\prime }}}}\wedge {\frac{{t^{\prime }}
}{{t}}}\right) ^{M} \frac {{\left( {{s}\vee {s^{\prime }}}\right) ^{4M}}}{{
\left( \left( {{s}\vee {s^{\prime }}}\right) ^{2}{+}\left\vert {z}
\right\vert ^{2}+\left\vert u\right\vert \right) {^{n+1+2M}}}}\frac{\left( {
{t}\vee {t^{\prime }}}\right) ^{2M}}{\left( {{t}\vee {t^{\prime }}+}
\left\vert {v}\right\vert \right) {^{1+2M}}}.
\end{eqnarray*}
Thus, we obtain that there exists a constant $C=C(M)$ depending only on $M$ such that
if $\left( {s}\vee {s^{\prime }}\right) ^{2}\leq {t}\vee {t^{\prime }},$
then
\begin{eqnarray*}
&&\left\vert \psi _{s,t}\ast _{\mathbb{H}^{n}}\ast K\ast_{\mathbb{H}^{n}}\psi _{s^{\prime },t^{\prime
}}(z,u)\right\vert \\
&\leq &C\left( {\frac{{s}}{{s^{\prime }}}}\wedge {\frac{{
s^{\prime }}}{{s}}}\right) ^{2M}\left( {\frac{{t}}{{t^{\prime }}}}\wedge {
\frac{{t^{\prime }}}{{t}}}\right) ^{M}  {\frac{\left( {{s}\vee {s^{\prime }}}\right) ^{2M}}{\left( {{s}\vee
{s^{\prime }}+}\left\vert {z}\right\vert \right) {^{2n+2M}}}}{\frac{\left( {{
t}\vee {t^{\prime }}}\right) ^{M}}{\left( {{t}\vee {t^{\prime }}+}\left\vert
{u}\right\vert \right) {^{1+M}}}},
\end{eqnarray*}
and if $\left( {s}\vee {s^{\prime }}\right) ^{2}\geq {t}\vee {t^{\prime }},$
then
\begin{eqnarray*}
&&\left\vert \psi _{s,t}\ast _{\mathbb{H}^{n}}\ast K\ast_{\mathbb{H}^{n}} \psi _{s^{\prime },t^{\prime }}(z,u)\right\vert\\
&\leq &C\left( {\frac{{s}}{{s^{\prime }}}}\wedge {\frac{{s^{\prime }}}{{s}}}
\right) ^{M}\left( {\frac{{t}}{{t^{\prime }}}}\wedge {\frac{{t^{\prime }}}{{t
}}}\right) ^{M} {\frac{\left( {{s}\vee {s^{\prime }}}\right) ^{M}}{\left( {{s}\vee {
s^{\prime }}+}\left\vert {z}\right\vert \right) {^{2n+M}}}}{\frac{\left( {{s}
\vee {s^{\prime }}}\right) ^{M}}{\left( {{s}\vee {s^{\prime }}+}\sqrt{
\left\vert {u}\right\vert }\right) {^{2+2M}}}}.
\end{eqnarray*}
These estimates imply that if $M>\alpha_1\vee \alpha_2,$
\begin{eqnarray*}
&&s^{-\alpha_1}t^{-\alpha_2}|\psi_{s,t}\ast T(f)(z,u)|\\
&\leq &C s^{-\alpha_1}t^{-\alpha_2}\int_0^\infty\int_0^\infty \left( {\frac{{s}}{{s^{\prime }}}}\wedge {\frac{{s^{\prime }}}{{s}}}
\right) ^{M}\left( {\frac{{t}}{{t^{\prime }}}}\wedge {\frac{{t^{\prime }}}{{t
}}}\right) ^{M}(s')^{\alpha_1}(t')^{\alpha_2}\frac{ds'}{s'}\frac{dt'}{t'}\|f\|_{Lip_{flag}}\\
&\leq& C \|f\|_{Lip_{flag}},
\end{eqnarray*}
which yields the proof of the claim.

We now extend $T$ to $Lip_{flag}^\alpha$ as follows. First, if $f\in Lip_{flag}^\alpha$ then, as mentioned above, there exists a sequence $\{f_n\}_{n\in\mathbb Z}\in L^2\cap Lip_{flag}^\alpha$ such that $f_n$ converges to $f$ in the distribution sense and $||f_n||_{Lip_{flag}^\alpha}\leq C||f||_{Lip_{flag}^\alpha}.$ It follows from the claim that $$||T(f_n)-T(f_m)||_{Lip_{flag}^\alpha}\leq C||f_n-f_m||_{Lip_{flag}^\alpha}$$ and hence $T(f_n)$ converges in the distribution sense. We now define $$T(f)=\lim\limits_{n\rightarrow\infty} T(f_n)$$ in the distribution sense. We obtain, by Theorem \ref{r1} and the above claim,
\begin{eqnarray*}||T(f)||_{Lip_{flag}^\alpha}&\lesssim &\sup\limits_{s,t>0,(z,u)\in\mathbb H^n}s^{-\alpha_1}t^{-\alpha_2}|\psi_{s,t}\ast T( f)(z,u)|
\\&\lesssim &\sup\limits_{s,t>0,(z,u)\in\mathbb H^n}s^{-\alpha_1}t^{-\alpha_2}|\lim\limits_{n\rightarrow\infty} \psi_{s,t}\ast T( f_n)(z,u)|
\\&\lesssim &\liminf\limits_{n\rightarrow\infty}\sup\limits_{s,t>0,(z,u)\in\mathbb H^n}s^{-\alpha_1}t^{-\alpha_2}| \psi_{s,t}\ast T( f_n)(z,u)|
\\&\lesssim &\liminf\limits_{n\rightarrow\infty}\|f_n\|_{Lip_{flag}^\alpha}
\\&\lesssim &\|f\|_{Lip_{flag}^\alpha}.\end{eqnarray*}
The proof of Theorem \ref{r2} is concluded.

\end{document}